\documentclass[11pt,reqno]{amsart}
\usepackage{amsmath,amssymb,amsthm,mathrsfs}
\usepackage[margin=1.1in]{geometry}
\usepackage{hyperref}

\usepackage{tikz}
\usepackage{enumitem}
\usetikzlibrary{arrows.meta}
\usetikzlibrary{decorations.pathreplacing}

\newtheorem{theorem}{Theorem}[section]

\newtheorem{corollary}[theorem]{Corollary}
\theoremstyle{definition}

\newtheorem{question}[theorem]{Question}

\newcommand{\R}{\mathbb{R}}
\newcommand{\C}{\mathbb{C}}

\newcommand{\Vol}{\operatorname{Vol}}
\newcommand{\Area}{\operatorname{Area}}

\newcommand{\Pin}{\Pi}

\title[A dichotomy for minimal submanifolds]{A dichotomy for minimal submanifolds}
\thanks{The authors were partially supported by NSF DMS Grants 2405393 and
2304684.}

\author{Tobias Holck Colding}
\author{William P.\ Minicozzi II}
\address{Department of Mathematics, Massachusetts Institute of Technology,
77 Massachusetts Avenue, Cambridge, MA 02139}
\email{colding@math.mit.edu and minicozz@math.mit.edu}

\begin{document}

\begin{abstract}
We survey a circle of recent results showing that a minimal submanifold obeys a
dichotomy: either it \emph{fills up space}, spreading out like a
space--filling curve, or it is \emph{confined}, and confinement forces 
quantitative restrictions.  On the confined side these restrictions are both
geometric and function--theoretic: Euclidean volume growth, an optimal rate of
convergence of the density, complex--curve rigidity for stable surfaces in
$\R^4$, and a Liouville theorem forcing slowly growing harmonic functions on
minimal disks to be constant.  A prototypical restriction is Euclidean volume growth, which in these results is \emph{forced by the geometry} rather than assumed.  
The mechanism is a volume doubling theorem that converts geometric confinement into quantitative rigidity: a stationary integral varifold
trapped in a thin slab at a given scale cannot double its volume by more than a
universal factor.  We explain how this principle, and the height--excess bounds
behind it, produce Euclidean volume growth for submanifolds of sublinearly
growing height in every dimension and codimension, the optimal density rate in a
slab, the complex--curve and Liouville rigidity above, a higher--codimension
Bernstein theorem for disks, one--sided volume bounds, and an optimal stable
Bernstein theorem in all dimensions generalizing Moser, Bombieri--De
Giorgi--Miranda, Caffarelli--Nirenberg--Spruck and Ecker--Huisken.

We  describe three regimes with sharp borderlines. In narrow cones the
volume grows polynomially with exponent $n+O(\varepsilon^2)$, and the order
$\varepsilon^2$ is attained; in wide cones the growth can be arbitrarily fast.
For holomorphic curves in $\mathbb{C}^2$, there are also three regimes and
 the borderline is determined by
the K\"ahler angle $\theta$ of the plane, with sharp thresholds
$\tan(\theta/2)$ and $\cot(\theta/2)$.
\end{abstract}

\maketitle

\section{From global rigidity to local estimates}

A minimal surface is a critical point of the area functional, or equivalently a surface whose mean curvature vanishes identically. Physically, minimal surfaces model soap films spanning wire frames, while mathematically they provide some of the simplest nonlinear geometric objects. Their study lies at the intersection of the calculus of variations, geometric analysis, and geometric measure theory, and has long driven the development of each of these subjects. Beginning with the pioneering work of Euler and Lagrange in the eighteenth century, which gave birth to the calculus of variations, the theory of minimal surfaces has time and again generated ideas and techniques that later became fundamental throughout regularity theory, geometric measure theory, differential geometry, and the analysis of nonlinear elliptic partial differential equations.

One of the central themes in elliptic PDEs is that rigidity often has a hidden origin. 
Strong conclusions are initially obtained under seemingly indispensable hypotheses, only to be discovered 
later to follow from substantially weaker assumptions. Much of the history of elliptic regularity can be viewed as a search for the geometric or analytic mechanism behind rigidity.

The classical example comes from the theory of harmonic functions.  
Liouville's theorem asserts that a bounded harmonic function on all of
Euclidean space is constant, and more generally that harmonic functions with
sufficiently slow growth must be affine or polynomial. These are striking global rigidity theorems, but they only apply to solutions defined on all of Euclidean space.

The real breakthrough came when the focus shifted from proving global rigidity to understanding local behaviour. Interior gradient estimates, Harnack inequalities, and Schauder estimates apply on arbitrary balls and make no reference to the behaviour of the solution at infinity. In practice, these local estimates almost always imply the corresponding Liouville theorems, while simultaneously providing much more detailed information about arbitrary solutions. The global rigidity theorems remain true, but they are no longer the main story; they become consequences of a more powerful local theory.

The same philosophy reappeared in a different guise in the resolution of Hilbert's nineteenth problem by De Giorgi and Nash. Their work established the regularity of weak solutions of uniformly elliptic equations and transformed modern regularity theory.   A central ingredient in De Giorgi's proof is the control of a scale--invariant measure of how far the graph of a solution deviates from an affine plane, commonly known as De Giorgi's height--excess.   This geometric quantity measures the height of the graph above its best approximating plane.

The subsequent work of Almgren and Allard extended this philosophy from graphs to arbitrary sets satisfying elliptic equations in a weak sense. Their theory replaced graphs by stationary varifolds. In contrast to De Giorgi's approach, however, the regularity theory of stationary varifolds is formulated in terms of the tilt--excess, measuring the deviation of tangent planes from a fixed plane, rather than the weaker height--excess.

This distinction is fundamental. The tilt--excess is a substantially stronger quantity: it implicitly contains a scale--invariant volume bound, and essentially all of the classical regularity theory relies upon such a bound. By contrast, the height--excess merely measures how closely a set lies inside a thin slab and contains no obvious information about its volume.

From this perspective, the natural question is whether the passage from height--excess to tilt--excess is genuinely unavoidable. Does one really need to assume measure--theoretic information in order to obtain regularity, or can it itself be recovered from weaker geometric information?   
For minimal submanifolds, the stronger measure-theoretic information can be recovered from weaker geometric information.

The central ingredient is the volume doubling theorem, which shows that confinement inside a sufficiently thin slab automatically generates the missing volume bound: it converts the weak geometric information encoded by the height--excess into the measure--theoretic information required by the classical regularity theory.  The new mechanism may be summarized as

\[
\text{height--excess}
\Longrightarrow
\text{volume bound}
\Longrightarrow
\text{tilt--excess}.
\]

Thus the role played by De Giorgi's height-excess in the theory of graphs extends, after the volume bound, to arbitrary stationary varifolds and hence the stronger tilt--excess is no longer the fundamental hypothesis. The weaker geometric information encoded by the height--excess already generates the measure--theoretic information required by the classical regularity theory.

This observation has several consequences. Most notably, it shows that Euclidean volume growth is produced by the geometry rather than assumed.   This leads to the optimal rate of convergence of the density, Liouville theorems in arbitrary codimension, rigidity theorems for stable minimal submanifolds, and higher-dimensional analogues of phenomena previously understood only in three dimensions.

The remainder of this article develops this philosophy by showing that
minimal submanifolds naturally divide into two sharply contrasting classes.
Either they spread throughout the ambient space with essentially unrestricted
complexity, or they remain confined, and confinement forces strong
quantitative rigidity. The volume doubling theorem is the mechanism that connects these two worlds.

\section{The dichotomy}	\label{s:dichotomy}

Taken together, the results surveyed here point to a remarkably simple organizing principle for minimal submanifolds.   We first state it informally.

\begin{quote}
\textbf{Dichotomy.} A minimal submanifold either \emph{fills up space}---in the
sense that it wanders all over $\R^{n+k}$, much as a space--filling curve
does---or it is \emph{confined in space}, in which case it must satisfy strong
quantitative restrictions.
\end{quote}

The point is that confinement is expensive.  Minimality is a maximally
spread--out condition---the mean curvature vanishes, so the submanifold cannot
bend toward any preferred direction---and so if one nonetheless \emph{traps} it
in a bounded region of the ambient space, one pays for this with rigidity.  Even
in high dimension and high codimension, where minimal submanifolds are otherwise
essentially unrestricted, confinement produces sharp structure.

This dichotomy is already implicit in several classical results in three dimensions, and much of what follows can be viewed as a unification and extension of these ideas:

\begin{itemize}
\item \emph{The half--space theorem} (Hoffman--Meeks \cite{HoffmanMeeks}, \cite{CM-distance}): a
      properly immersed nonplanar minimal surface in $\R^3$ cannot lie in a
      half--space.  A minimal surface that is confined to one side of a plane is
      forced to be a plane.
\item \emph{Embedded minimal disks in $\R^3$} (Colding--Minicozzi
      \cite{CM-disks1,CM-disks2,CM-disks3,CM-disks4}): a sequence of embedded
      minimal disks that does not stay confined near a point degenerates into a
      foliation by parallel planes away from a singular axis, with the disks spiraling like a  
      helicoid; those that do stay confined are graphical and controlled.
\item \emph{The one--sided curvature estimate} (Colding--Minicozzi
      \cite{CM-disks4}): an embedded minimal disk in $\R^3$ that lies on one side
      of a plane, close to it, is uniformly graphical---confinement to one side of
      a plane forces curvature bounds.
\end{itemize}

In this survey we describe how the same dichotomy operates in \emph{all}
dimensions and codimensions.  The unifying restriction that confinement forces is
\textbf{Euclidean volume growth}, and the unifying mechanism is a
\textbf{volume doubling theorem} \cite{CM-confined}.  We then survey the
consequences: control of density, the equivalence of the natural notions of
closeness to a plane, an optimal rate of convergence of the density in a slab,
complex--curve rigidity in $\R^4$, one--sided bounds, a Liouville theorem for
minimal disks \cite{CM-liouville}, and Bernstein--type theorems.  In the opposite
direction, the examples of Colding--Mart\'in--Minicozzi \cite{CMM} show that
without confinement the area can grow arbitrarily fast, so the confinement
hypotheses cannot be dropped.

Throughout, the single most important quantity is the \emph{scale--invariant
density}
\[
   \Theta(x,r)\;=\;\frac{\Vol\!\big(B_r(x)\cap\Sigma\big)}{\omega_n\, r^n},
   \qquad \Sigma^n\subset \R^{n+k},
\]
where $\omega_n$ is the volume of the unit ball in $\R^n$.  By the monotonicity formula $\Theta (x,r)$ is monotone nondecreasing in $r$ for a minimal
$\Sigma$; see \cite{Simon,ColdingMinicozziBook}.  Density controls both the local
regularity theory and the global structure of a minimal submanifold.   A bound
on the density is almost always assumed---often implicitly.  A central message here
is that such a bound need not be assumed: for confined minimal submanifolds it is
\emph{implied by the geometry}.

\section{The volume doubling theorem}

To make this dichotomy quantitative one needs a mechanism converting confinement into volume control.  The following volume doubling theorem converts geometric confinement into volume control and underlies many of the results discussed below.  

\begin{theorem}[Volume doubling, Colding--Minicozzi \cite{CM-confined}]
\label{thm:doubling}
There exist $\delta_0>0$ and $C$ depending only on $n$ so that if
$\Sigma^n\subset\R^{n+k}$ is a proper stationary integral varifold and
$B_{4\,r}\cap\Sigma$ is contained in a slab of height $\delta_0\, r$, then
\[
   \Vol\!\big(B_{2\,r}\cap\Sigma\big)\;\le\; C\,\Vol\!\big(B_r\cap\Sigma\big).
\]
\end{theorem}

The theorem says that geometric confinement forces quantitative control of volume.   More precisely, a stationary integral varifold confined inside a sufficiently thin slab cannot more than double its volume by a universal factor. It assumes only geometric confinement--closely related to De Giorgi's height--excess--and produces the scale--invariant volume bound that had previously been assumed rather than proved.

Two features are essential:
\begin{enumerate}
\item The doubling constant $C$ is \emph{universal} (it depends only on $n$), not
      on the submanifold, and in particular is independent of the codimension $k$.
\item Nothing is assumed about the density a priori; the volume bound is produced,
      not imposed.  The slab hypothesis may vary from scale to scale, which
      permits fractal behavior.
\end{enumerate}

\subsection{Iteration and Euclidean volume growth}
Once confinement has been shown to produce volume control, the first natural question is whether it produces the correct asymptotic growth.  When the slab condition holds at \emph{all} large scales, one may iterate.  This
is exactly the case of \emph{sublinearly growing height}: after a rigid motion
there are $\alpha<1$ and $r_0>0$ with
\[
   \Sigma\setminus B_{r_0}\subset
   \Big\{\,x_{n+1}^2+\cdots+x_{n+k}^2<\big(x_1^2+\cdots+x_n^2\big)^{\alpha}\,\Big\}.
\]
Iteration alone yields only a polynomial bound on the volume growth. With additional ideas and a more refined analysis, one obtains Euclidean volume growth: 

\begin{theorem}[Euclidean volume growth, Colding--Minicozzi \cite{CM-confined}]
\label{thm:evg}
A complete proper stationary integral varifold of any dimension and codimension
whose height grows sublinearly has Euclidean volume growth: there is $C$ with
\[
   \Vol\!\big(B_r\cap\Sigma\big)\;\le\; C\,r^n\qquad\text{for all }r\ge1.
\]
Moreover the constant $C$ is determined by the volume of a single sufficiently
large ball, of radius depending only on $r_0,\alpha,n$.
\end{theorem}

The exponent is exactly $n$: the confined minimal submanifold has Euclidean
volume growth, and this is the sharp quantitative restriction promised by the
dichotomy.  The generalized Enneper surfaces and examples of Kapouleas, \cite{Kapouleas}, show that
the multiplicative constant can be arbitrarily large.  
The sublinear hypothesis cannot be weakened to linear height: already in  
$\mathbb{R}^3$ there are stationary integral varifolds in a cone of aperture
$\varepsilon$ whose volume grows at least like $r^{2+\frac45\varepsilon^2}$
\cite{CM-height}; see Section~\ref{s:S11new}.

What does Euclidean volume growth buy us? Before answering that question, let us first see why it is the right growth condition.

\subsection{Why Euclidean volume growth is the right invariant}
There are several ways to see that Euclidean volume growth is not merely a technical convenience; it is the natural global condition governing the theory.   A useful analogy: for a smooth
$f:\R^n\to\R$ with negative gradient flow
$X_x(0)=x,\ \tfrac{d}{dt}X_x(t)=-\nabla f(X_x(t))$,
the most basic qualitative question is whether $f$ is bounded below; boundedness
below prevents the flow from running off to $-\infty$ and organizes the dynamics.
In the minimal--surface setting, Euclidean volume growth plays exactly this role:
it is the global condition that rules out the degenerate, space--filling behavior
and makes structure possible.  Two classical instances make this concrete.

\begin{itemize}
\item \textbf{Surfaces in $\R^3$.} For a minimal surface in $\R^3$ of finite
      genus, Euclidean area growth is equivalent to \emph{finite total
      curvature}, and such surfaces have been completely understood for decades
      (Jorge-Meeks \cite{JorgeMeeks}, Osserman \cite{Osserman}, Schoen \cite{Schoen}, and others).
\item \textbf{Complex submanifolds and algebraic geometry.} Wirtinger
      \cite{Wirtinger} observed in 1936 that every complex submanifold of $\C^N$
      is area minimizing.  Over more than fifteen years in the 1950s--60s, using value
      distribution theory (a complex analogue of Hadamard's three--circles
      theorem) together with work of Hermann Weyl, the combined efforts of
      Rutishauser, Stoll and Bishop \cite{Rutishauser,Stoll,Bishop} established that a complex
      submanifold of $\C^N$ with Euclidean volume growth is \emph{algebraic},
      i.e.\ the zero set of finitely many polynomials.  Euclidean volume growth is
      exactly the dividing line between transcendental wandering and algebraic
      rigidity.
      This raises the question of which geometric confinements force Euclidean
growth, and hence algebraicity; for cones the answer is sharp and depends on 
the K\"ahler angle of the plane (Theorem~\ref{t:11p3}). 
\end{itemize}

Once Euclidean volume growth is available, the next natural question is how it changes our notion of closeness to a plane.

\section{From geometric to measure-theoretic closeness}

Because density plays a central role in both the local and global theory, it is worth isolating what the doubling
theorem buys us at the level of closeness to a plane.  Consider three natural
notions that a set $\Sigma$ is close, on $B_r$, to a ball in an $n$--plane:
\begin{enumerate}
\item \textbf{Hausdorff close:} $\Sigma\cap B_r$ is Hausdorff--close to a disk in
      the plane;
\item \textbf{Measure--theoretically close:} $\Sigma\cap B_r$ is close in measure
      to a disk in the plane;
\item \textbf{Varifold close:} $\Sigma$ together with its tangent planes is close,
      as a varifold, to the plane with multiplicity one (roughly a $W^{1,2}$
      condition).
\end{enumerate}
For a stationary varifold \emph{with a fixed volume bound}, these three notions are
all equivalent \cite{Allard,Simon}.  The natural question is what happens without
an a priori volume bound.

\begin{question}
Without assuming a volume bound, are these notions still comparable?  Does the
weakest of them---Hausdorff closeness---control the others?
\end{question}

The answer, which is the conceptual content of Theorem~\ref{thm:doubling}, is
yes: Hausdorff closeness places $\Sigma$ in a thin slab, doubling then produces
the missing volume bound, and once one has a volume bound the classical
equivalences \cite{Allard,Simon} apply.  Thus the weakest and most flexible
hypothesis turns out to imply them all.

\subsection{Height--excess versus tilt--excess}
This distinction has a long history in elliptic PDE and geometric measure theory.
De Giorgi \cite{DeGiorgi} and Nash \cite{Nash} independently resolved Hilbert's
nineteenth problem in the late 1950s, establishing H\"older continuity of
derivatives of solutions of nonlinear elliptic equations.  A key ingredient in De
Giorgi's approach was control of a scale--invariant $L^2$ deviation of the graph
of a function from an affine plane---the \emph{height--excess}---and its decay at
smaller scales.  Almgren \cite{Almgren} and Allard \cite{Allard}, \cite{Delellis}, recast this in a
coordinate--free, varifold setting valid for general sets solving an elliptic
equation weakly, using instead a stronger, tangent--plane based quantity---the
\emph{tilt--excess}---which for graphs corresponds roughly to a $W^{1,2}$ bound.

The tilt--excess implicitly carries a scale--invariant \emph{volume bound}, and all existing arguments in this circle rely on such a
bound---yet many applications only supply control closer in spirit to De Giorgi's
original \emph{height}--excess.  The doubling theorem removes the gap: for confined
minimal submanifolds the weaker height--excess bound already implies a volume
bound, and hence a tilt--excess bound.  This is the technical source of the
applications that follow.

\section{The picture in three dimensions}

\subsection{Building blocks}
The volume doubling theorem is the mechanism underlying the higher-dimensional theory. In three dimensions, however, an even more detailed picture is available. There the same dichotomy is expressed not through volume doubling but through the structure theory of embedded minimal disks, built on the one-sided curvature estimate and reflected globally in the half-space theorem.  In three dimensions the
structure theory of Colding--Minicozzi
\cite{CM-disks1,CM-disks2,CM-disks3,CM-disks4,CM-genus}, \cite{Perez} shows that, locally, there are
exactly three essential building blocks:
\begin{enumerate}
\item planes and, more generally, graphs;
\item double spirals (pieces of the helicoid);
\item catenoid necks.
\end{enumerate}
A soap film is a stable minimal surface found in nature. 
Its geometry consists precisely of the three basic pieces discussed here: a flat or graphical sheet, 
one half of the double--spiral staircase of a helicoid, and not too large portion of the neck of a catenoid.

\subsection{The one--sided curvature estimate}
The technical heart of the $3$--dimensional theory is the \emph{one--sided
curvature estimate} \cite{CM-disks4}: an embedded minimal disk in a ball that lies
on one side of a plane and comes close to that plane is graphical over the plane
with uniformly bounded curvature on a definite smaller ball.

\begin{theorem}[One--sided curvature estimate, Colding--Minicozzi \cite{CM-disks4}]
There are $\varepsilon_0>0$ and $C$ so that the following holds.  If
$\Sigma^2\subset B_{2\,r}\cap\{x_3>0\}\subset\R^3$ is an embedded minimal disk with
$\partial\Sigma\subset\partial B_{2r}$, then each component of $B_r\cap\Sigma$
that comes within $\varepsilon_0\, r$ of the plane $\{x_3=0\}$ is a graph over that
plane with
\[
   \sup_{B_r\cap\Sigma}|A|^2 \;\le\; \frac{C}{r^2}.
\]
\end{theorem}

This is a direct manifestation of the dichotomy: a minimal disk that is
\emph{confined to one side of a plane} is forced to be graphical and flat.  If it
were not confined---if it wandered to both sides---it would spiral like a helicoid
and have large curvature.

\subsection{The half--space theorem}
At the global level, the Hoffman--Meeks half--space theorem \cite{HoffmanMeeks} captures the first manifestation of this picture.

\begin{theorem}[Half--space theorem \cite{HoffmanMeeks}]
A connected, properly immersed minimal surface in $\R^3$ contained in a
half--space is a plane.
\end{theorem}

The original proof of the half--space theorem, due to Hoffman and Meeks, \cite{HoffmanMeeks}, combines the strong maximum principle with a sliding catenoid argument. More recently, a conceptually different proof was found, \cite{CM-distance}. It is based on the distance to the minimal surface from points in the boundary of the half--space, together with the observation that suitably chosen functions of this distance are viscosity subsolutions. The argument is then completed by exploiting the parabolicity of the flat two-dimensional plane.

Again, confinement to a half--space is so restrictive for a minimal surface that only the flat ones survive. This striking rigidity is special to dimension three.  
In higher dimensions, the half--space theorem is known to fail: higher-dimensional catenoids lie in slabs, and many other minimal hypersurfaces are likewise   slab-contained. Holomorphic curves in $\mathbb{C}^2$ are an exception in higher codimension: by
Sibony--Wong \cite{SW}, a nonconstant proper holomorphic curve in a real half-space is a
complex line. This rigidity disappears in $\mathbb{C}^3$, and in
$\mathbb{C}^2$ it reappears quantitatively for cones (Section~\ref{s:S11new}).
The results described below may be regarded as quantitative, volume-theoretic refinements of this dichotomy.

The structure theory for embedded minimal disks in $\R^3$ developed in \cite{CM-disks1,CM-disks2,CM-disks3,CM-disks4,CM-genus} is not merely a local analogue of the half--space theorem, but a fundamentally stronger theory. Whereas the half--space theorem is a global rigidity result that applies only to complete minimal surfaces extending to infinity, the one--sided curvature estimate applies to an arbitrary local piece of a minimal surface and leads to a detailed geometric description of its structure.  At a high level, this is analogous to the distinction between a Liouville theorem for entire solutions and a local gradient estimate. In the absence of a sufficiently strong compactness theorem, estimates of the latter type are dramatically more powerful. In fact, the one--sided curvature estimate and the structure theorems of \cite{CM-disks1,CM-disks2,CM-disks3,CM-disks4,CM-genus} provide substantially more information than a local version of the half--space theorem could hope to provide.

\section{Without constraints, chaos: rapid area growth}

Even for surfaces \emph{without spatial constraints} minimal surfaces can behave
wildly.  The examples of Colding--Mart\'in--Minicozzi \cite{CMM} make this
precise, and show that the confinement hypotheses in the theorems above cannot be
removed.

\begin{theorem}[Colding--Mart\'in--Minicozzi \cite{CMM}]
The graph $\Sigma=\{(z,f(z)):z\in\C\}\subset\C^2=\R^4$ of the holomorphic
function $f(z)=\sin(e^z)$ is a smooth, properly embedded, \emph{area--minimizing}
(in particular stable) minimal surface, and there are $c,r_0>0$ with
\[
   \mathrm{Area}\big(B_r\cap\Sigma\big)\;\ge\; e^{c\,r}\qquad (r\ge r_0).
\]
Taking $f(z)=\sin(e^{z^2})$ gives instead $\mathrm{Area}(B_r\cap\Sigma)\ge
e^{c\,r^2}$.
\end{theorem}

\begin{theorem}[Colding--Mart\'in--Minicozzi \cite{CMM}]
There is a complete \emph{proper} conformal minimal immersion $X:\Omega\to\R^3$
of a simply connected surface with
\[
   \int_{X^{-1}(B_R)} dA_X\;\ge\; C\, e^{c\,R}
\]
for all large $R$; the radii can be chosen so that the growth is instead
Gaussian, $\ge C\, e^{c\,R^2}$.
\end{theorem}

These examples---even in $\R^3$ and $\R^4$---spread
out and fill up space: their density ratios blow up, so they violate volume
doubling on every large scale and are not confined to slabs across scales. Confinement to a cone of linear aperture does not prevent this: sums of
scaled catenoids, and proper holomorphic disks in wide cones, have arbitrarily
fast growth while staying in a cone \cite{CM-height}. 
 This
is the ``fills up space'' side of the dichotomy.  Everything below is about the
other side.

\subsection{Examples: how dimension changes confinement}
The sublinear--height hypothesis is satisfied by many natural minimal
submanifolds, and it is instructive that raising the dimension tends to make
confinement \emph{easier}, not harder.

\begin{itemize}
\item \textbf{Catenoids.} In $\R^3$ the height of the catenoid grows at the logarithmic rate, so it is
      not contained in any half--space, though it does have Euclidean area growth.  
      In higher dimensions the situation is
      strikingly different---the height of the higher--dimensional catenoid is
      \emph{bounded} \cite{Blair}, so it lies in a slab.  This is a first sign
      that the half--space theorem must fail for $n>2$: there really are
      minimal hypersurfaces trapped in a slab, and they fall squarely
      within the scope of Theorems~\ref{thm:doubling} and \ref{thm:evg}.
\item \textbf{Enneper surfaces.} Enneper's surface in $\R^3$ has polynomial height
      growth of order $\tfrac23$---sublinear, hence confined in the sense of
      Theorem~\ref{thm:evg}, though not contained in any slab.  Its
      higher--dimensional generalizations are conjectured to lie within a slab
      \cite{Choe}, which would again place them among the confined examples with
      bounded height rather than merely sublinear height.
\end{itemize}

Thus the passage from $n=2$ to $n>2$ converts the two borderline
$3$--dimensional examples---the catenoid and Enneper's surface---from
``sublinearly confined'' to (conjecturally, for Enneper) ``slab confined,''
exactly the regime where the sharper slab results of the next sections apply.

\section{An optimal rate of convergence in a slab}

When the submanifold lies in a fixed slab, one gets not only Euclidean volume
growth but an optimal rate at which the density approaches its limit, together with
integrality of that limit.

\begin{theorem}[Colding--Minicozzi \cite{CM-confined}]	\label{t:confined}
There exist $C$ and $R_0$ depending only on $n$ so that if $\Sigma^n\subset
\{\,x_{n+1}^2+\cdots+x_{n+k}^2\le 1\,\}$ is a complete proper stationary integral
varifold, then there is an integer $N$ so that for all $r\ge R_0$ and all $p$ in
the axis $\{x_{n+1}=\cdots=x_{n+k}=0\}$,
\[
   \big(1-C\, r^{-2}\big)\,N\;\le\;\frac{\Vol\big(B_r(p)\cap\Sigma\big)}
   {\Vol\big(B_r\subset\R^n\big)}\;\le\; N.
\]
\end{theorem}

The $r^{-2}$ rate is sharp, as the example of a slightly tilted $n$--plane in the
slab shows, and the integrality of $N$ follows from Allard's integrality theorem
once the density and excess are controlled.  The proof rests on a weighted
integrability of the tilt--excess, $\sum_j\int_\Sigma|\nabla^T x_{n+j}|^2\,
|x|^{2-n}<\infty$, which improves the naive logarithmic bound that follows once one has shown Euclidean volume growth.

One of the first consequences of Euclidean volume growth is function-theoretic rigidity.

\section{A Liouville theorem and a Bernstein theorem for disks}

Euclidean (here quadratic) area growth is exactly the condition under which a
function--theoretic Liouville phenomenon appears.  For minimal \emph{disks} in any
codimension this takes the following sharp form \cite{CM-liouville}.

\begin{theorem}[Liouville theorem, Colding--Minicozzi \cite{CM-liouville}]
\label{thm:liouville}
Let $\Sigma^2\subset\R^N$ be a properly immersed minimal surface that is
topologically a disk, with quadratic area growth $\Area (B_r\cap\Sigma)\le C_a\, r^2$.  If
$u$ is a harmonic function on $\Sigma$ whose negative part grows slowly,
\[
   -C\big(1+|x|^{\alpha}\big)\le u(x)
   \quad\text{for some constant $C$ and some }
   \alpha<\frac{-\log\big(1-e^{-24\,C_a}\big)}{\log 2},
\]
then $u$ is constant.
\end{theorem}

Since the coordinate functions are harmonic, one must have $\alpha<1$ for any such
statement, and the admissible range above is indeed $<1$.  Combining
Theorem~\ref{thm:liouville} with the volume growth estimates of
\cite{CM-confined} gives a Bernstein theorem for minimal disks of slowly growing
height in any codimension.

\begin{corollary}[Bernstein theorem for disks, Colding--Minicozzi
\cite{CM-liouville}]
Given $N$, there is $\alpha>0$ so that if $\Sigma$ is a properly immersed minimal
disk in $\R^N$ contained in
$\{x\in\R^N:\sum_{i=3}^N|x_i|\le C(|x_1|^\alpha+|x_2|^\alpha+1)\}$ for some
constant $C$, then $\Sigma$ is the $x_1-x_2$ plane.
\end{corollary}

The result is optimal, and the three classical surfaces calibrate it exactly:
\begin{itemize}
\item \emph{Enneper's surface} has quadratic area growth with $C_a=3\,\pi$, tangent
      cone at infinity a plane of multiplicity three, and a nonconstant harmonic
      coordinate growing at rate $\tfrac23$; higher--order Enneper surfaces make
      $\alpha\to0$ as $C_a\to\infty$.  So some sub--linear restriction is genuinely
      needed.
\item \emph{The catenoid} shows that the topological (disk) hypothesis cannot be
      dropped even in $\R^3$: its logarithmically growing vertical coordinate is a
      nonconstant harmonic function.
\item \emph{The helicoid} has cubic area growth and a logarithmically growing
      harmonic function, showing the quadratic area bound is necessary.
\end{itemize}

The same oscillation estimate yields an effective, quantitative version.

\begin{theorem}[Uniform H\"older regularity, Colding--Minicozzi
\cite{CM-liouville}]
If $\Sigma\subset B_{2\,r}\subset\R^N$ is a compact immersed minimal disk with
$\partial\Sigma\subset\partial B_{2\,r}$ and $\Area(\Sigma)\leq 4\,C_a\,r^2$, then there is $\alpha=\alpha(C_a)>0$ so that any harmonic $u$ on $\Sigma$ satisfies for $x$, $y$ in the same component of $B_s\cap\Sigma$ with $0<s<r$, 
\[
   |u(x)-u(y)|\le C\,\|u\|_{L^1(\Sigma)}\,\Big(\frac{s}{r}\Big)^\alpha.
\]
\end{theorem}

A point of interest is that the Gauss--Bonnet theorem is never used, so the
techniques are available in higher dimensions as well.

\section{Stable minimal submanifolds}

Stable minimal surfaces---those that locally minimize area to second order---are
the soap films one sees in nature: minimal graphs, part (but not all) of the
catenoid, half of the helicoid, half of Enneper's surface.  For this class the
dichotomy takes an especially clean form.   

The strongest rigidity appears for stable minimal submanifolds.

\subsection{An optimal stable Bernstein theorem in all dimensions}

\begin{theorem}[Stable Bernstein, Colding--Minicozzi \cite{CM-confined}]	\label{t:SBT}
Let $\Sigma^n\subset\R^{n+1}$ be a complete, properly immersed, two--sided stable
minimal hypersurface.  If $\Sigma$ has sublinearly growing height, then $\Sigma$
is a hyperplane.
\end{theorem}

The proof uses the volume bound of Theorem~\ref{thm:evg} to make the
scale--invariant tilt--excess small enough to apply Bellettini's $\varepsilon$--
regularity \cite{Bellettini,SchoenSimon}, and then a reverse Poincar\'e argument.
The sublinear growth is necessary: stable minimal cones, such as the Simons cone
in $\R^8$, have exactly linear height growth. (Linear height is also the threshold for volume growth alone; see Section~\ref{s:S11new}.)  
 Theorem \ref{t:SBT} generalizes the classical
Bernstein theorems for minimal graphs, where stability and Euclidean volume
growth hold automatically:
\begin{itemize}
\item Moser (1961) \cite{Moser}: bounded $|\nabla u|$;
\item Bombieri--De Giorgi--Miranda (1969) \cite{BDGM}: sublinear growth of $u$;
\item Trudinger (1972) \cite{Trudinger}: a new proof of Moser's theorem;
\item Caffarelli--Nirenberg--Spruck (1988) \cite{CNS}: a generalization of
      Moser's theorem;
\item Ecker--Huisken (1990) \cite{EckerHuisken}: sublinear growth of $\nabla u$.
\end{itemize}
There is also a version with compact boundary, giving an asymptotic description
(disjoint graphs converging to constants) and the improved curvature decay
$|A|^2\le C\,|x|^{2\,\alpha-4}$ outside a compact set \cite{CM-confined}.

\subsection{Stable surfaces in $\R^4$ are complex curves}
In higher codimension, Micallef \cite{Micallef} proved that a complete oriented
parabolic stable minimal surface in $\R^4$ is holomorphic for some orthogonal
complex structure.  Combined with the volume bounds above, this gives:

\begin{theorem}[Colding--Minicozzi \cite{CM-confined}]	\label{t:cmMM}
Let $\Sigma^2\subset\R^4$ be an oriented stable stationary integral varifold with
at most finitely many singular points.  If $\Sigma$ is contained in a sublinearly
growing tubular neighborhood of some two--plane, then it is a complex curve for
some orthogonal complex structure on $\R^4$ and indeed algebraic.
\end{theorem}

The sublinear growth gives quadratic area growth (Theorem~\ref{thm:evg}), hence
parabolicity by the logarithmic cutoff argument, so Micallef's theorem applies.
Many holomorphic curves satisfy the hypotheses, e.g.\ $z_1^p=z_2^q$ with $p\ne q$
in $\C^2=\R^4$.

Under a parabolicity assumption,
 the sublinear hypothesis can be relaxed to
confinement in a cone of aperture $\varepsilon<1$; see Corollary~\ref{c:old113}.

\section{One--sided volume bounds}

The results above assume a two--sided (slab or tube) bound.  Confinement to one
side of a hyperplane still yields a doubling, but---as tilted planes show---one
must decrease the height of the slab as the radius grows.

\begin{theorem}[Colding--Minicozzi \cite{CM-confined}]
There are $C_1,C_2,C_3>0$ depending on $n$ and $k$ so that if $z_0>0$ satisfies
$r>C_1 z_0$ and $C_{4r}\cap\Sigma\subset\{x_{n+j}>0,\ j=1,\dots,k\}$, then
\[
   \Vol\big(C_{2r,\,C_2 z_0}\cap\Sigma\big)\;\le\;
   C_3\,\Vol\big(C_{r,\,z_0}\cap\Sigma\big),
\]
where $C_{r,b}=\{x_1^2+\cdots+x_n^2<r^2\}\cap\{x_{n+j}<b\}$.
\end{theorem}

Thus one gains volume doubling as one goes further out, at the necessary cost of
looking closer to the bounding hyperplane.

\section{Three regimes with sharp thresholds}	\label{s:S11new}

Theorem~\ref{thm:evg} assumes the height grows sublinearly. The first case it does not
cover is linear height, that is, confinement to a thickened cone with aperture $\varepsilon$
around an $n$-plane
$\Pi$:
\[
S_{\varepsilon,K}(\Pi)=\{x : |P^\perp x|\le \varepsilon |Px| + K\},
\qquad C_\varepsilon(\Pi)=S_{\varepsilon,0}(\Pi),
\]
where $P,P^\perp$ are the orthogonal projections onto $\Pi,\Pi^\perp$.
 In \cite{CM-height} it is shown that linear height is  the sharp
 borderline.

\begin{figure}[ht]
\centering
\begin{tikzpicture}[
   x=1.15cm,
   y=1.15cm,
   scale=.7,
   every node/.style={font=\small},
   line cap=round,
   line join=round
]
   \coordinate (O) at (0,0);
   \def\xmax{3.45}
   \def\slope{0.42}
   \def\offset{0.68}

   \fill[black!10]
      (-\xmax,{\slope*\xmax+\offset})
      --(0,\offset)
      --(\xmax,{\slope*\xmax+\offset})
      --(\xmax,{-\slope*\xmax-\offset})
      --(0,-\offset)
      --(-\xmax,{-\slope*\xmax-\offset})
      --cycle;

   \draw[->,thick] (-3.72,0)--(3.78,0)
      node[below right] {$\Pin$};
   \draw[->,thick] (0,-2.30)--(0,2.36)
      node[above] {$\Pin^\perp$};

   \draw[thick]
      (-\xmax,{\slope*\xmax+\offset})
      --(0,\offset)
      --(\xmax,{\slope*\xmax+\offset});
   \draw[thick]
      (-\xmax,{-\slope*\xmax-\offset})
      --(0,-\offset)
      --(\xmax,{-\slope*\xmax-\offset});

   \draw[densely dashed,black!65]
      (-\xmax,{\slope*\xmax})--(O)--(\xmax,{\slope*\xmax});
   \draw[densely dashed,black!65]
      (-\xmax,{-\slope*\xmax})--(O)--(\xmax,{-\slope*\xmax});

   \draw[<->,thin] (.25,.04)--(.25,{\offset-.04})
      node[midway,right=2pt] {$K$};

   \node at (2.25,1.58) {$S_{\varepsilon,K}(\Pin)$};
   \node at (2.45,.43) {$C_\varepsilon(\Pin)$};
   \fill (O) circle (1.2pt);
\end{tikzpicture}
\caption{The dashed lines bound the cone
$C_\varepsilon(\Pin)$ with aperture $\varepsilon$.  The cone is contained in $S_{\varepsilon,K}(\Pin)$ where the constant $K$ permits a bounded
transverse error.}
\label{f:bounded-error-cone}
\end{figure}

\subsection{Stationary varifolds}

\begin{theorem}[Colding--Minicozzi \cite{CM-height}]	\label{t:11p1}
Let $\Sigma^n\subset\mathbb{R}^{n+k}$ be a stationary integral varifold and
$V(r)=\mathrm{Vol}(B_r\cap\Sigma)$.
\begin{enumerate}
\item If the height of $\Sigma$ over $\Pi$ grows sublinearly, then $V$ has
Euclidean growth (Theorem~\ref{thm:evg}).
\item There are $\delta_n,\delta_n'>0$ so that if
$\Sigma\subset S_{\varepsilon,K}(\Pi)$ with $\varepsilon<\delta_n$, then $V$
grows polynomially; if $\varepsilon<\delta_n'$, then
$V(r)\le c\,(1+r)^{\,n+c_n\varepsilon^2}$.
\item For $\varepsilon>\delta_n$ there are examples with super-polynomial
growth, and $\delta_n\le \max_{t\ge0} t/\cosh t\approx 0.66$.
\end{enumerate}
\end{theorem}

\begin{figure}[ht]
\centering
\begin{tikzpicture}[
    x=1.76cm,
    y=1.15cm,
    scale=.82,
    every node/.style={font=\small},
    line cap=round,
    line join=round
]
    \coordinate (O) at (0,0);

    \fill[black!12]
       (O)--(1.4,.56)--plot[domain=1.4:5.25,samples=100]
          (\x,{.56+.3394*(sqrt(\x-1.22)-sqrt(.18))})--(5.25,0)--(O);
    \shade[bottom color=black!45,top color=white]
       (O)--(5.25,2.75)--(5.25,4.70)--(0,4.70)--cycle;

    \draw[->,thick]
       (-.85,0)--(5.28,0);
    \node[below right] at (4.85,-.03) {$\Pi$};
    \draw[->,thick]
       (0,-1.05)--(0,4.78) node[above] {$\Pi^\perp$};

    \draw[thick,black!55]
       (O)--(1.4,.56)
       plot[domain=1.4:5.25,samples=100]
          (\x,{.56+.3394*(sqrt(\x-1.22)-sqrt(.18))});
    \draw[thick,black!75] (O)--(5.25,2.75);
    \node[rotate=20,above=2pt] at (5.55,2.8) {$\delta_n$};
    \fill (O) circle (1.3pt);

    \node[align=center] at (3.7,.42)
       {Euclidean growth};
    \node[align=center,rotate=12,font=\footnotesize] at (3.7,1.50)
       {Polynomial growth};
    \node[align=center] at (2.95,3.05)
       {Super polynomial growth};

    \draw[decorate,
          decoration={brace,amplitude=5pt,mirror},
          thick]
       (5.43,.04)--(5.43,1.05)
       node[midway,right=8pt,align=left]
       {sublinear\\[-1pt] height};

    \draw[decorate,
          decoration={brace,amplitude=5pt,mirror},
          thick]
       (5.43,1.18)--(5.43,2.68)
       node[midway,right=8pt,align=left]
       {narrow\\[-1pt] cones};
\end{tikzpicture}
\caption{Three regimes of volume-growth for \emph{stationary varifolds}.}
\label{f:real-growth-regimes}
\end{figure}

\begin{figure}[ht]
\centering
\includegraphics[width=0.5\textwidth]{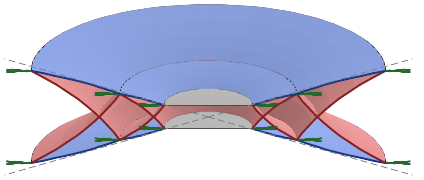}
\caption{Half of the stationary varifold in Theorem \ref{t:11p2}, cut along a vertical plane
through the $x_3$-axis.  It consists of two flat disks (gray) and pieces of
catenoids (blue and red) that meet along circles and stay in the cone
$|x_3|\le\varepsilon r$ (dashed); the catenoidal ends are only indicated, in
green, on the cut.}
\label{f:varifold3d}
\end{figure}

The exponent in (2) is sharp, even in codimension one:

\begin{theorem}[Colding--Minicozzi \cite{CM-height}]	\label{t:11p2}
For every $\varepsilon\in(0,1/10)$ there is a stationary integral 2-varifold
in $S_{\varepsilon,1}(\Pi)\subset\mathbb{R}^3$, invariant under rotations and
consisting of two flat disks and pieces of catenoids meeting along circles,
with $V(r)\ge r^{2+\frac45\varepsilon^2}$ for $r$ large.
\end{theorem}

The mechanism is a balancing argument along the junction circles at radii
$2^k$. Along a catenoid the vertical flux $r\sin\vartheta$ is constant, so
each sheet arrives at the next circle flatter than it left. A flatter sheet
pulls harder horizontally, so more sheets must leave a circle than arrive at
it. The number of sheets therefore grows by a factor  on each dyadic scale, which yields the exponent
$2+\tfrac45\varepsilon^2$. Compare Theorem~\ref{t:confined}: in a slab (bounded height)
the density converges to an integer at the optimal rate, whereas in a cone it
can diverge.  

\subsection{Holomorphic curves}

For holomorphic curves the borderlines are linear. Identify
$\mathbb{R}^4=\mathbb{C}^2$. The K\"ahler angle $\theta\in[0,\pi/2]$ of a
two-plane $\Pi$ is given by $\cos\theta=|\omega(e_1,e_2)|$ where $\omega$ is the K\"ahler form and $e_1 , e_2$ is
an orthonormal frame for $\Pi$.
Thus $\theta=0$
for a complex line and $\theta=\pi/2$ for a Lagrangian plane.

\begin{theorem}[Holomorphic thresholds \cite{CM-height}]	\label{t:11p3}
For nonconstant proper holomorphic curves in $\mathbb{C}^2$ (see Figure \ref{fig:kahler-angle-regions}):
\begin{enumerate}
\item if $\varepsilon<\tan(\theta/2)$, none lies in $S_{\varepsilon,K}(\Pi)$
for any $K<\infty$;
\item if $\tan(\theta/2)\le\varepsilon<\cot(\theta/2)$, such curves exist,
and every one is affine algebraic with quadratic area growth and parabolic
domain;
\item if $\varepsilon>\cot(\theta/2)$, there are complete proper holomorphic
immersions of the disk into $C_\varepsilon(\Pi)$ with arbitrarily fast area
growth.
\end{enumerate}
\end{theorem}

\begin{figure}[ht]
\centering
\begin{tikzpicture}[
    line width=0.6pt,
    >={Stealth[length=2.6mm,width=2mm]},
    every node/.style={font=\small}
  ]

  \def\angHi{35}      
  \def\angLo{19}      
  \def\Lhi{7.4}       
  \def\Llo{6.6}       

  \fill[gray!12] (0,0) -- (\angHi:\Lhi) -- (0, {\Lhi*sin(\angHi)}) -- cycle;      
  \fill[gray!25] (0,0) -- (\angHi:\Lhi) -- (\angLo:\Llo) -- cycle;             
  \fill[gray!6]  (0,0) -- (\angLo:\Llo) -- ({\Llo*cos(\angLo)},0) -- cycle;      
  
  \draw[->] (0,0) -- (8.0,0) node[pos=0.88, below=5pt] {$\Pi$};
  \draw[->] (0,0) -- (0,5.4)  node[pos=0.92, left=5pt]  {$\Pi^{\perp}$};
  
  \draw (0,0) -- (\angHi:\Lhi) node[right, xshift=2pt] {$\cot\dfrac{\theta}{2}$};
  \draw (0,0) -- (\angLo:\Llo) node[right, xshift=2pt] {$\tan\dfrac{\theta}{2}$};

  \path (0,0) -- (\angHi:\Lhi)
        node[sloped, above, pos=0.52] {Rapid growth examples};
  \path (0,0) -- (\angLo:\Llo)
        node[sloped, above, pos=0.70] {Algebraic, quadratic growth};
  \path (0,0) -- (6.2,0)
        node[above, pos=0.62] {Nonexistence};

\end{tikzpicture}
\caption{The three \emph{holomorphic} aperture regimes.  The thresholds are determined by the K\"ahler angle $\theta$ of the two-plane $\Pi$.}
\label{fig:kahler-angle-regions}
\end{figure}

\subsection{Stable surfaces in $\mathbb{R}^4$}

Combining (2) with Micallef's theorem \cite{Micallef} extends Theorem~\ref{t:cmMM} from sublinear to
linear height, at the cost of assuming parabolicity:

\begin{corollary}[\cite{CM-height}]	\label{c:old113}
A complete oriented properly immersed stable minimal surface in
$\mathbb{R}^4$ of parabolic type that lies in $S_{\varepsilon,K}(\Pi)$ with
$\varepsilon<1$ is an affine algebraic curve for some orthogonal complex
structure $J$, has quadratic area growth, and satisfies
$\theta_J(\Pi)\le 2\arctan\varepsilon$.
\end{corollary}

In particular, a stable counterexample to quadratic growth in narrow cones
must be nonparabolic, and so has area growth faster than $R^2\log R$.

\section{Sketch of the proof of doubling}

Finally, we indicate the proof of Theorem~\ref{thm:doubling} when  
$\Sigma^n\subset\R^{n+1}$, $n>2$, is a smooth minimal hypersurface, contained in a
narrow slab $\{|x_{n+1}|\le \delta R\}$ with $\delta>0$ small.  The argument illustrates how confinement (thinness of the slab) is converted into a
volume bound.

\smallskip
\noindent\textbf{Step 1: the cylindrical distance.}
Let $s>0$ be the distance to the $x_{n+1}$--axis, \[s^2 = x_1^2+\cdots+x_n^2\, ,\] and
set $V(r)=\Vol(\{s\le r\}\cap\Sigma)\approx \Vol(B_r\cap\Sigma)$.  We want
$V(2\,r)\le C\,V(r)$.

\smallskip
\noindent\textbf{Step 2: tilt--excess as flatness.}
The quantity $|\nabla^T x_{n+1}|^2$ is the squared deviation of $T_x\Sigma$ from
$\{x_{n+1}=0\}$; its average is the tilt--excess.

\smallskip
\noindent\textbf{Step 3: comparing $s$ with a harmonic model.}
On $\Sigma$ one has
\[
   -(n-1)\,(n-2)\,|\nabla^T x_{n+1}|^2
   \;\le\; s^{n}\,\Delta_\Sigma\,\!\big(s^{2-n}\big)
   \;\le\; (n-2)\,|\nabla^T x_{n+1}|^2 ,
\]
and integrating this bounds $\int_{\Sigma_{2R}}|\nabla^T s|^2$ by
\[n^2\, 2^{2n-1}V(R)+(n-1)\,2^{n-1}\int_{\Sigma_{R,2R}}|\nabla^T x_{n+1}|^2\,  .\]

\smallskip
\noindent\textbf{Step 4: the normal derivative is small.}
Since $\Sigma$ is nearly horizontal, $|\nabla^\perp s|^2 \le |\nabla^T x_{n+1}|^2$.

\smallskip
\noindent\textbf{Step 5: a subharmonic barrier.}
Set $h = s^{2-n} + 2^{\,n-1}\, n\,(n-2)\,x_{n+1}^2/r^n$.  If the slab is thin enough,
then on $\{\tfrac r2\le s\le 4\,r\}$,
\[
   \varepsilon\, s^{2-n}\ge|s^{2-n}-h|\,  ,\]
\[   \Delta_\Sigma\, h \ge 2\,n\,(n-2)\,\frac{|\nabla^T x_{n+1}|^2}{r^n}.
\]

\smallskip
\noindent\textbf{Step 6: bounding the tilt--excess by $V(r)$.}
Using $h$ as a barrier and integrating $\Delta_\Sigma\, h$ bounds
$\int_{\{r\le s\le 2\,r\}}|\nabla^T x_{n+1}|^2$ in terms of $V(r)$.  Thinness of the
slab is used decisively here to control the boundary terms.

\smallskip
\noindent\textbf{Step 7: Pythagoras and conclusion.}
Since $s$ has unit gradient, $1=|\nabla^T s|^2+|\nabla^\perp s|^2$, so
$V(2\,r)=\int_{\Sigma_{2\,r}}(|\nabla^T s|^2+|\nabla^\perp s|^2)$.  By Steps 3--6 both
terms are bounded by the tilt--excess, which is bounded by $V(r)$.  Hence
$V(2\,r)\le C\,V(r)$. \qed

\section{Concluding remarks}

Section~\ref{s:S11new} shows that the dichotomy of Section~\ref{s:dichotomy} is the coarse
form of a finer picture.  As the confinement is relaxed, a minimal
submanifold passes through three regimes: a rigid regime, an intermediate
regime in which quantitative control survives, and a flexible regime in which
that control is lost.  The regimes are separated by sharp thresholds, and the
same pattern appears for stationary varifolds and for holomorphic curves.

For stationary varifolds, the confinement is measured by the growth of the
height over an $n$-plane $\Pi$ (Theorem~\ref{t:11p1}).  Sublinear height
forces Euclidean volume growth.  In a narrow cone of aperture $\varepsilon$,
the growth is polynomial with exponent $n+O(\varepsilon^2)$, and the order
$\varepsilon^2$ is attained (Theorem~\ref{t:11p2}).  Beyond a threshold
$\delta_n$, the growth can be super-polynomial, and in cones of aperture larger
than $\max_{t\ge0}t/\cosh t\approx0.66$,
 sums of scaled catenoids grow
arbitrarily fast.  The first threshold is sharp in the strongest sense:
Euclidean growth holds whenever the height is sublinear, and there are
examples without it as soon as the height is linear, however small the slope.

For holomorphic curves in $\C^2$, both thresholds are explicit and are
determined by the K\"ahler angle $\theta$ of $\Pi$
(Theorem~\ref{t:11p3}).  In cones of aperture below $\tan(\theta/2)$, there
are no curves at all; between $\tan(\theta/2)$ and $\cot(\theta/2)$, every
curve is affine algebraic with quadratic area growth; above $\cot(\theta/2)$,
the area can grow arbitrarily fast.  The lower threshold is exactly the
aperture at which the cone first contains a complex line.  The two extreme
planes show how much the thresholds encode.  If $\Pi$ is a complex line, then
$\theta=0$ and the algebraic regime fills every aperture.  If $\Pi$ is
Lagrangian, then $\theta=\pi/2$, both thresholds equal $1$, and the algebraic
regime disappears.

The two pictures have the same shape, but the regimes sit differently.  For
varifolds, there is no empty regime, since $\Pi$ itself is stationary and lies
in every cone around it; the rigid regime is Euclidean growth, and it ends as
soon as the height becomes linear.  For holomorphic curves, the rigid regime is
nonexistence, and Euclidean growth occupies a whole interval of linear
apertures, reflecting the additional rigidity of complex curves.

Identifying the thresholds is more than a matter of completeness, because the
thresholds carry information.  Their values record the structure that governs
the problem, and the examples at a threshold show how rigidity is lost.
Corollary~\ref{c:old113} illustrates the first point: both its hypothesis and
its conclusion can be read off from the holomorphic thresholds.  Since
$\cot(\theta/2)\ge1$ for every plane, an aperture $\varepsilon<1$ stays below
the upper threshold for whichever orthogonal complex structure $J$ makes the
surface holomorphic, so the surface lies in the algebraic regime.  The bound
$\theta_J(\Pi)\le2\arctan\varepsilon$ says precisely that $\varepsilon$ lies
above the lower threshold $\tan(\theta_J(\Pi)/2)$.  The aperture of the cone
thus restricts which complex structures can occur.  For varifolds, the
examples of Theorem~\ref{t:11p2} show how Euclidean growth is lost at the
first threshold: balancing along the junction circles forces the number of
catenoidal sheets to grow from one dyadic scale to the next.  The same linear
threshold appears in the stable Bernstein theorem (Theorem~\ref{t:SBT}), where
the obstruction is instead the Simons cone.

We expect this to be part of a larger picture.  In many variational problems,
relaxing a constraint moves the solutions from rigidity, through a regime of
controlled structure, to flexibility, and the thresholds between the regimes
carry essential information.  Classical examples suggest how widespread
this may be.  Harmonic functions on $\R^n$, the critical points of the
Dirichlet energy, are constant if they grow sublinearly, form
finite-dimensional spaces of polynomials if they grow polynomially, and form
infinite-dimensional spaces if they may grow exponentially.  The first
threshold is again sublinear growth, sharp because of the coordinate
functions, and both rigid statements persist on manifolds of nonnegative Ricci
curvature, by the gradient estimate of Cheng--Yau \cite{CY} and the finite
dimensionality theorem of \cite{CM-harmonic}.

\end{document}